\documentclass[conference]{IEEEtran}
\IEEEoverridecommandlockouts
\usepackage{cite}
\usepackage{amsmath,amssymb,amsfonts}
\usepackage{algorithmic}
\usepackage{svg}
\usepackage{graphicx}
\usepackage{textcomp}
\usepackage{xcolor}
\usepackage{orcidlink}

\hypersetup{
    colorlinks,
    linkcolor={red!50!black},
    citecolor={blue!50!black},
    urlcolor={blue!80!black}
}

\makeatletter
\newcommand{\linebreakand}{%
  \end{@IEEEauthorhalign}
  \hfill\mbox{}\par
  \mbox{}\hfill\begin{@IEEEauthorhalign}
}
\makeatother

\def\BibTeX{{\rm B\kern-.05em{\sc i\kern-.025em b}\kern-.08em
    T\kern-.1667em\lower.7ex\hbox{E}\kern-.125emX}}

\begin{document}

\title{Escape from an Orbiting Pursuer with a Nonzero Capture Radius\\
\thanks{This research was funded by Air Force Research Laboratory under FA8650-20-2-5853 and 21RQCOR084. 
Distribution Statement A: Approved for public release.  Distribution is unlimited.  Case number AFRL-2023-4713.
}
}

\author{\IEEEauthorblockN{Braulio Mora}
\IEEEauthorblockA{\textit{Mechanical and Aerospace Engineering Department} \\
\textit{The University of Texas at Arlington}\\
Arlington, TX, USA \\
braulio.mora@mavs.uta.edu} \and
\IEEEauthorblockN{Alexander Von Moll \orcidlink{0000-0002-7661-5752}}
\IEEEauthorblockA{\textit{Control Science Center} \\
\textit{Air Force Research Laboratory}\\
Wright-Patterson AFB, OH USA \\
alexander.von\_moll@afrl.af.mil} \and 

\IEEEauthorblockN{Isaac Weintraub \orcidlink{0000-0003-2209-2500}}
\IEEEauthorblockA{\textit{Control Science Center} \\
\textit{Air Force Research Laboratory}\\
Wright-Patterson AFB, OH, USA \\
isaac.weintraub.1@afrl.af.mil} \linebreakand

\IEEEauthorblockN{David W. Casbeer}
\IEEEauthorblockA{\textit{Control Science Center} \\
\textit{Air Force Research Laboratory}\\
Wright-Patterson AFB, OH USA \\
david.casbeer@afrl.af.mil} \and  

\IEEEauthorblockN{Animesh Chakravarthy}
\IEEEauthorblockA{\textit{Mechanical and Aerospace Engineering Department } \\
\textit{The University of Texas at Arlington}\\
Arlington, TX, USA \\
animesh.chakravarthy@uta.edu}
}

\maketitle

\begin{abstract}
This paper explores a  multi-agent containment problem, where a fast evader, modeled having constant speed and using constant heading, attempts to escape a circular containment region that is orbited by a slower pursuer with a nonzero capture radius.
The pursuer is constrained to move along the edge of the containment region and seeks to capture the evader.
This paper presents an in-depth analysis of this pursuer-evader containment scenario.
First, multiple types of capture conditions for a single-pursuer case are analyzed defining the worst-case initial position for the pursuer.
Second, a parametric study is performed to demonstrate the effects of speed ratio, capture radius, and initial location of the evader.
Finally, a reachability analysis is performed to investigate the viable escape headings and reachable regions by the evader.
This work provides a foundation for the analysis of escape under more general evader inputs as well as a multiple-pursuer version of the scenario.
\end{abstract}


%%%%%%%%%%%%%%%%%%%%%%%%%%%%%%%%%%%%%%%%%%%%%%%
\section{Introduction}
\label{sec:intro}
The field of differential games and optimal control has proven to be an effective tool in the analysis of adversarial scenarios with applications in vehicle motion planning~\cite{mirzaei2021optimal}, military combat~\cite{ardema1985combat, isaacs1965differential}, medicine~\cite{buratto2015hiv}, and economics~\cite{friedman2013differential}.
The pursuit-evasion game is one of many fundamental problems that can be related to real world scenarios~\cite{basar1982dynamic}.
A well formulated introduction to pursuit-evasion problems formulated as differential games is made by~\cite{isaacs1965differential}.
Extensions to the basic setup of the problem have been made and added to the literature (c.f., \cite{weintraub2020introduction}).
The Tag-Chase game, the homicidal chauffeur, and the lady in the lake are some examples of scenarios studied over the years.

An extension to the ``Lady in the Lake" problem is being considered in this paper.
The Lady in the Lake problem originally posed by Martin Gardner in a mathematics column~\cite{gardner2008lady} (later analyzed in detail in~\cite{breakwell1977zero-sum,basar1982dynamic}) is an instance of a reach-avoid problem, which is simply a special type of pursuit-evasion problem where the lady (evader) attempts to escape the lake (a containment region) without being captured by a pursuer standing at the edge of the lake.
In the scenario considered here, a fast evader ($E$) is located within a circular containment region and moves with a constant speed and fixed heading angle.
A slower pursuer ($P$) orbits the containment region and is constrained to clockwise (CW), or counter-clockwise (CCW) motion along the edge of the region.
It is endowed with a finite (nonzero) capture radius.
Capture is said to occur if $E$ comes within the capture radius of $P$.
The goal of $E$ is to escape the containment region while $P$ attempts to capture $E$.
Note that, since it is assumed that the Evader is faster, escape (under general heading control inputs) against a single pursuer is always possible.
This study is mainly concerned with straight-line motion of the evader as such behavior often arises as an extremal control in such scenarios as this~\cite{basar1982dynamic}.

A similar problem setup is analyzed in~\cite{yan2017escape-avoid}, in which an evader seeks to escape a circular region orbited by multiple pursuers in an evenly distributed formation.
There, the evader is assumed to be slower than the pursuers, so there are regions of the state space in which the latter can prevent the former from escaping.
Additionally, that work considers capture to occur if the evader and any pursuer become collocated (i.e., point capture).
Here, the evader is given the advantage of speed while the pursuer is given an advantage in having a nonzero capture radius.
As will be shown throughout the remainder, the inclusion of the nonzero capture radius significantly complicates the analysis of the scenario.

Several other works have analyzed scenarios in which an agent (or team of agents) wish to contain an adversary within an environment (or prevent its entry into a protected area).
For example, the works~\cite{garcia2021cooperative,chakravarthy2020cooperative} seek to develop cooperation among a team of agents to prevent the escape of an evader.
Like this paper, the works~\cite{yan2017defense,garcia2020optimal} involve a pursuit-evasion game which takes place inside a circle in which the evader seeks to escape prior to capture; there, the pursuer(s) was not constrained to the circle boundary.
The classic Lion and Man differential game~\cite{flynn1973lion} also takes place in a circular play region, but the game (typically) does not end if the evader is able to reach the boundary of the play region.
Containment of an adversary is also an objective of the herding scenario presented and analyzed in~\cite{chipade2019herding}.
In~\cite{shishika2018local-game} a team of fast defenders constrained on a circular region are tasked with preventing intruders from reaching the circle from outside (by capture via point capture).
In~\cite{fu2023defending} a slower defender with a nonzero capture radius attempts to either capture an intruder or prevent it from entering the protected area indefinitely.
Finally, in~\cite{rivera2023pursuer} as a nonlinear control synthesis method was developed to enable a team of slower pursuers to prevent an evader from reaching its target.

The contributions of the paper are as follows:
\emph{i)} identification of all possible capture configurations,
\emph{ii)} computation of the worst-case pursuer initial condition for a given evader heading (i.e., the one in which the pursuer starts the furthest away),
\emph{iii)} a parametric study of the capture configurations and worst-case pursuer initial conditions,
\emph{iv)} a reachability analysis from the perspective of the evader which identifies safe headings for the evader to take.
The utility of these results is demonstrated on the analysis of a two-pursuer variant of the problem.

The remainder of the paper is organized as follows.
Section~\ref{sec:problem_formulation} describes the mathematical formulation of the problem.
Section~\ref{sec:worst_case_capture_for_one_pursuer}, which represents the bulk of the study, contains the analysis of the various capture configurations as well as the worst-case pursuer initial position.
Sections~\ref{sec:parametric_study} and~\ref{sec:reachability_analysis} contain the parametric study and evader reachability analysis, respectively.
The paper is concluded in Section~\ref{sec:conclusion} with a few remarks on some extensions for which the results presented herein will be valuable.

\section{Problem Formulation}
\label{sec:problem_formulation}

\begin{figure}[ht]
    \centering
    \includegraphics[width=0.9\linewidth]{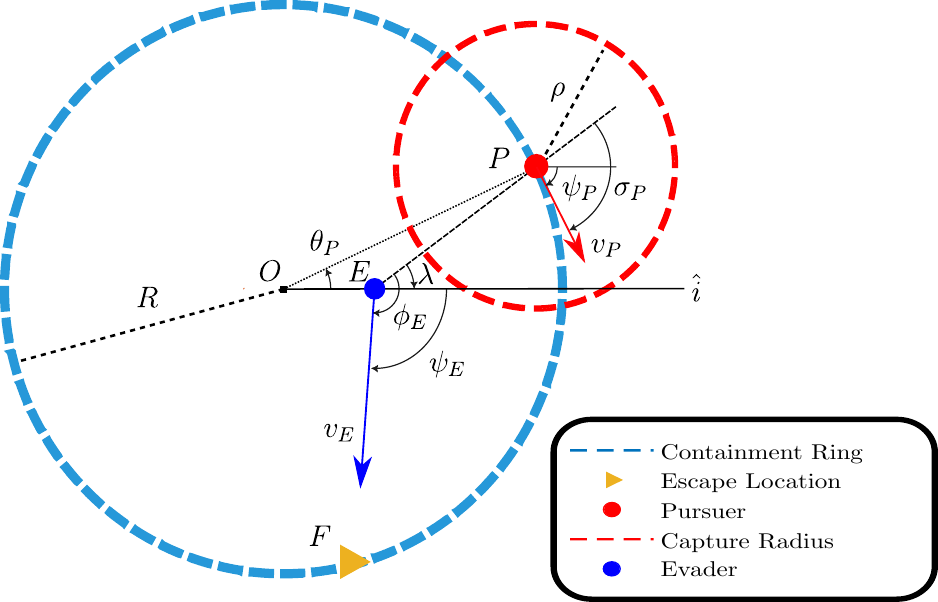}
    \caption{Basic Scenario with One Pursuer}
    \label{fig:basic_geo}
\end{figure}

The pursuer-evader scenario in two-dimensional space is illustrated in Fig. \ref{fig:basic_geo}. In this scenario, a fast evader, $E$, is located within a containment circular region of radius $R$ and attempts to escape. A slow pursuer, $P$, is located along the containment region edge, or ring, and is constrained to circular motion along the ring only. $P$ has a nonzero capture radius, $\rho > 0$, and tries to capture $E$ before it escapes. We consider $E$ and $P$ to have constant speeds $v_E$ and $v_P$, respectively. The speed ratio is defined as $\gamma=v_P/v_E<1$ and, without loss of generality, we normalize the speeds so that $v_E = 1$ and $v_P = \gamma$. 

In Cartesian coordinates, the dynamics of $P$ are
\begin{equation}
\begin{aligned}
\label{eq:dyn_P_cart}
    \dot{x}_P &= v_P\cos{\psi_P},         &   \dot{y}_P &= v_P\sin{\psi_P} .
    \end{aligned}
\end{equation}
where $(x_P,y_P)$ denotes the position of the pursuer, and $\psi_P$ denotes the heading angle of $P$. As aforementioned, $P$, moves in circular motion, thus, the heading angle of $P$ is
\begin{align}
\label{eq:p_heading}
    \psi_P = \theta_P + a\frac{\pi}{2}  .
\end{align}
Here, $\theta_P$ represents $P$'s angle with respect to the $\hat{i}$ axis, and  $a \in \{-1,1\}$ represents a quantity which depends on the direction of motion of $P$.  We use $a = -1$  for clockwise motion (CW) of $P$ and $a=1$ for counter-clockwise motion (CCW) motion of $P$. 

Similarly, the dynamics of $E$ in Cartesian coordinates are described as: 
\begin{align}
\label{eq:dyn_E_cart}
    \dot{x}_E &= v_E\cos{\psi_E},         &   \dot{y}_E &= v_E\sin{\psi_E} .
\end{align}
where $(x_E,y_E)$ denotes the position and $\psi_E$ denotes the heading angle of the evader. The heading angle of $E$ is assumed to be constant with time for this formulation. 

Given the circular geometry of the scenario, it is only natural to transform the Cartesian dynamics into relative polar coordinates. The equations of motion are developed using the geometry of Fig. \ref{fig:basic_geo}.  Defining the line of sight (LOS) distance between $E$ and $P$ as $R_{PE}$, and the LOS angle as $\lambda$, we have the following: %We begin by calculating these two states from the Cartesian locations:   
\begin{align}
\label{eq:R_lamda_eq}
    R_{PE} &= \sqrt{(x_E-x_P)^2+(y_E-y_P)^2},\\   
    \lambda &= \tan^{-1}{\frac{y_E-y_P}{x_E-x_P}}     .
\end{align}

The angles between the velocity of the evader and the LOS, $\phi_E$, %can be found by subtracting the line of sight angle from the heading of the evader. Similarly, the angle 
and that between the velocity of the pursuer and the LOS, $\sigma_P$, are as follows: 
\begin{align}
\label{eq:phi_sigma_eq}
    \phi_E &= \psi_E-\lambda,\\   
    \sigma_P &= \psi_P-\lambda    .
\end{align}
Substituting $\psi_P$ from \eqref{eq:p_heading}, we calculate the rates of change for the LOS angle, and the relative velocity angles as: 
\begin{alignat}{2}
\label{eq:angles_rates_eq}
    &\dot{\lambda}& &= \frac{v_e\sin{\phi_E}-v_P\sin{\sigma_P}}{R_{PE}},\\
    &\dot{\phi}_E& &= -\dot{\lambda},\\   
    &\dot{\sigma}_P& &= \dot{\theta}_P-\dot{\lambda},\\
    &\dot{\theta}_P& &= a\frac{v_P}{R}.
\end{alignat}
  
The distance rate is given by: 
\begin{align}
\label{eq:distance rate_eq}
    \dot{R}_{PE} &= v_E\cos{\phi_E}-v_p\cos{\sigma_P}.
\end{align}
Given the equations of motion shown above, we can discard the $\dot{\lambda}$ state since its motion is captured by $\dot{\phi}_E$ and $\dot{\sigma}_P$. A final relative polar space is described as: 
\begin{equation}
    \begin{bmatrix}
    \dot{R}_{PE}\\
    \dot{\theta}_P\\
    \dot{\phi}_E\\
    \dot{\sigma}_P
    \end{bmatrix} = 
    \begin{bmatrix}
    v_E\cos{\phi_E}-v_p\cos{\sigma_P}\\
    \frac{a}{R}v_P\\
    -\dot{\lambda}\\
    \dot{\theta}_P-\dot{\lambda}
    \end{bmatrix}
\end{equation}
Note that this is not a minimal polar representation. Some of the states are redundant, but we include them here as a continence for the purpose of the analysis. 

\section{Worst-case Capture for One Pursuer}
\label{sec:worst_case_capture_for_one_pursuer}

This section covers the general scenario for one pursuer, the different types of capture conditions and the analysis of each of them, as well as the final process to find the capture condition that allows for the maximum angular range for the pursuer's initial condition, i.e. the initial condition that is furthest away from $E$'s escape location along the containment ring. Initial conditions that correspond to this worst-case scenario require the pursuer to travel the largest distance before capturing the evader, that is, $\theta_{P0} - \theta_F$ is the largest.

\subsection{Capture Types}
The following subsection covers the multiple types of capture that are defined and analyzed. The goal is to obtain equations that give the most limiting (i.e. the worst-case scenario) initial location of the pursuer to achieve capture based on the speed ratio, $\gamma$, $P$'s capture radius, $\rho$, $E$'s initial location along the x-axis, $r$, and $E$'s aim point angle, $\theta_F$. 
Without loss of generality it is assumed in this section that the pursuer moves clockwise (that is, $a = -1$).
\subsubsection{Point Capture}
\begin{figure}
    \centering
    \includegraphics[width=0.45\linewidth]{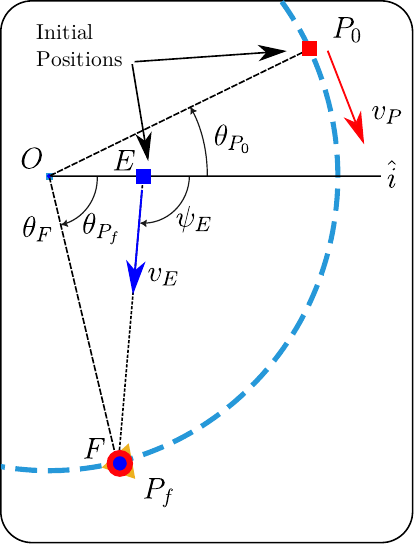}
    \caption{Point Capture Geometry}
    \label{fig:pc_geo}
\end{figure}

Let us first analyze the simplest case of point capture where the capture radius of $P$ is $\rho=0$. Point capture occurs when the position of $E$ and $P$ are collocated on top of each other. Since $P$ is constrained to motion along the ring, capture will occur when $E$ is positioned along the ring as shown in Fig. \ref{fig:pc_geo}. We can calculate the time to capture by dividing the distance travelled by $E$ or $P$ from its initial position to the exit point, $d$, by its velocity.  
\begin{equation}
\label{eq:time_eq}
t_c = d/v.
\end{equation}

Starting with the evader, and assuming the aim point angle, $\theta_F$, is known, we can obtain a value for the distance travelled, $d$, or according to Fig. \ref{fig:pc_geo}, segment $\overline{EF}$, using the Law of Cosines. Substituting for the respective values we obtain the following: 
\begin{equation}
\label{eq:EF}
\overline{EF}=\sqrt{R^2 + \overline{OE}^2 - 2R\overline{OE}\cos{\theta_F}}.
\end{equation}
We now substitute $\overline{EF}$ in \eqref{eq:time_eq} for $d$ to obtain the following: 
\begin{equation}
\label{eq:tc_E}
t_c = \overline{EF}/v_E.
\end{equation}

Similarly, we can obtain the time of capture of $P$ by using the arc length equation to solve for $d$ as follows: 
\begin{equation}
\label{eq:P1Pf}
d = R(\theta_{P0} - \theta_{Pf}),
\end{equation}
where $\theta_{P0}$ and $\theta_{Pf}$ are the initial and final angles of $P$, respectively. Just as before, we substitute \eqref{eq:P1Pf} into \eqref{eq:time_eq} appropriately to obtain:
\begin{equation}
\label{eq:tc_P}
t_c = R(\theta_{P0} - \theta_{Pf})/v_P.
\end{equation}

Now, we can equate \eqref{eq:tc_E} and \eqref{eq:tc_P} and solve for the initial angle of $P$, $\theta_{P0}$. This results in the following simple equation: 
\begin{equation}
\label{eq:IC_pc}
\theta_{P0} = \theta_{Pf} -\frac{\gamma\overline{EF}}{R},
\end{equation}
where the point capture condition, $\theta_{Pf} = \theta_F$, the speed ratio, $\gamma=v_P/v_E$, and $a\in \{-1,1\}$ depending on the direction of $P$.

Additionally, we can obtain the heading of $E$ again by using the Law of Cosines: 
\begin{equation}
\label{eq:psi_e}
\psi_E=\pi \mp \arccos{\frac{\overline{OE}^2+\overline{EF}^2-R^2}{2\overline{OE}\overline{EF}}}.
\end{equation}

\subsubsection{Exit Point Capture}
\begin{figure}
    \centering
    \includegraphics[width=0.45\linewidth]{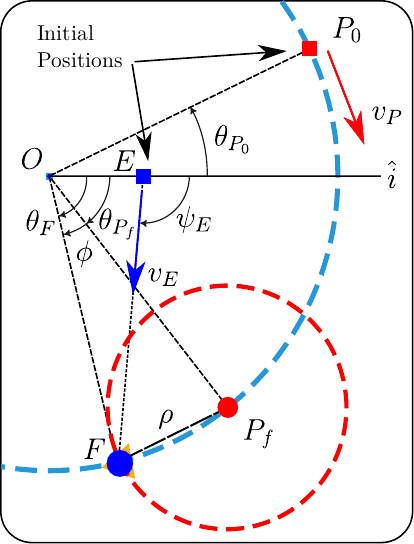}
    \caption{Exit Point Capture Geometry}
    \label{fig:epc_geo}
\end{figure}

\begin{figure}[ht]
    \centering
    \includegraphics[width=0.8\linewidth]{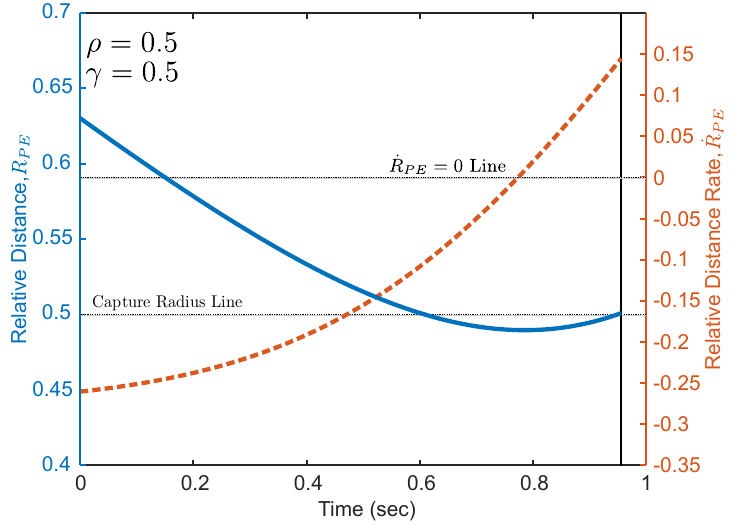}
    \caption{Exit Point Capture Relative Distance and Relative Distance Rate }
    \label{fig:epc_distance}
\end{figure}

For this condition, $P$ is given a nonzero capture radius, $\rho > 0$. Exit point capture (EXC) occurs when $P$ captures $E$ when $E$ is located along the containment ring as shown in Fig. \ref{fig:epc_geo}. The same methods as those used for point capture before can be used for this capture condition. The capture times of $E$ and $P$ can be set equal to each other to obtain: 
\begin{equation}
    \label{eq:tc_equal}
    \overline{EF}/v_E = R(\theta_{P0} - \theta_{Pf})/v_P.
\end{equation}
Now, $\theta_{Pf}$ is not simply $\theta_F$ but note that for this capture condition, there is an offset angle $\phi$ as shown in Fig. \ref{fig:epc_geo}. This offset is calculated using the chord formula as such: 
\begin{equation}
    \label{eq:phi_epc}
    \phi=2\arcsin{\frac{\rho}{2R}}.
\end{equation}
The final angle of $P$ is now: 
\begin{equation}
    \label{eq:theta_pf_epc}
    \theta_{Pf}=\theta_F - 2\arcsin{\frac{\rho}{2R}}.
\end{equation}
Now, we can calculate the initial position of $P$ for an EXC condition as follows: 
\begin{equation}
\label{eq:theta_p0_EXC}
\theta_{P0} = \theta_F - \frac{\gamma\overline{EF}}{R} - 2\arcsin{\frac{\rho}{2R}}.
\end{equation}
\emph{Analysis:} After analyzing the EXC condition on different scenarios, it was found that for some given combination of speed ratio, capture radius, and $E$'s heading angle, and due to the geometry of the problem, the evader gets captured before it reaches the exit point $F$ as shown in the blue curve in Fig. \ref{fig:epc_distance}. In this figure, the area below the capture radius line represents the evader going into the pursuer's capture zone. Note that the blue curve crosses this line at two points which indicates that there is a point earlier in time where $E$ could be captured, and an initial condition of $P$ further along the ring that could accomplish this. We can determine that this is not the worst-case solution for this scenario, and this is not always the case. A comparison is made later.

\subsubsection{Tangent Capture}
\begin{figure}
    \centering
    \includegraphics[width=0.45\linewidth]{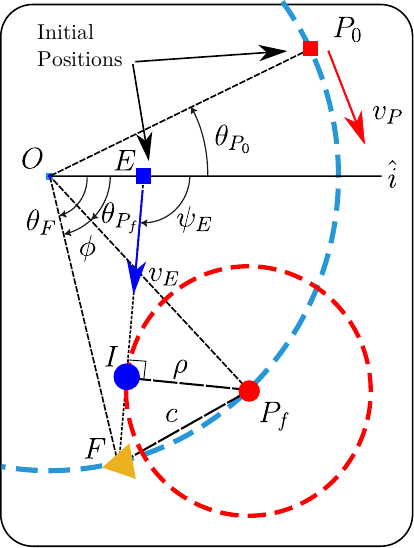}
    \caption{Tangent Capture Geometry}
    \label{fig:tan_geo}
\end{figure}

Tangent capture (TAC) occurs when the the pursuer intersects the $E$'s path at a tangent point. The vector loop closure method is used to find the distance travelled by the evader before capture. For this method, we treat the geometry shown in Fig. \ref{fig:tan_geo} as a four-bar mechanism consisting of segments $\overline{OE}, \overline{EI}, \overline{IP}$ and $\overline{PO}$. Note that segments $\overline{EI}$ and $\overline{IP}$ have a perpendicular relation, and that $\overline{IP}$ is equal to the capture radius, $\rho$, of the pursuer. As such, these are treated as a single bar. The vector loop closure consists of summing the $\hat{i}$ and $\hat{j}$ components of the vector loop such that the result is zero. 
\begin{align}
    \overline{OE}\hat{i} + \overline{EI}\hat{i} + \overline{IP}\hat{i} + \overline{PO}\hat{i} &= 0, \\
    \overline{OE}\hat{j} + \overline{EI}\hat{j} + \overline{IP}\hat{j} + \overline{PO}\hat{j} &= 0.
\end{align}
Substituting in the respective component values we obtain the following: 
\begin{align}
    x_E + \overline{EI}\cos{\psi_ E} + \rho\cos\left(\psi_ E + \pi/2\right) &= R\cos{\theta_P}, \\
    y_E + \overline{EI}\sin{\psi_ E} + \rho\sin\left(\psi_ E + \pi/2\right) &= R\sin{\theta_P}.
\end{align}
The above equations are simplified by taking into account the trigonometric identities $\cos\left(x+\pi/2\right)=-\sin{x}$ and $\sin\left(x+\pi/2\right)=\cos{x}$. Additionally, we assume $E$ to always lie in the $\hat{i}$ axis such that $y_E = 0$. 
\begin{align}
    x_E + \overline{EI}\cos{\psi_ E} - \rho\sin{\psi_ E} = R\cos{\theta_P}, \\
    \overline{EI}\sin{\psi_ E} + \rho\cos{\psi_ E} = R\sin{\theta_P}.
\end{align}
Both equations are now squared and summed. Using the trigonometric identity $\sin^2{x} + \cos^2{x} = 1$, we are able to cancel $\theta_P$ and obtain:  
% \begin{multline}
%     (x_E + \overline{EI}\cos{\psi_ E} - \rho\sin{\psi_ E})^2 \\ + (\overline{EI}\sin{\psi_ E} + \rho\cos{\psi_ E})^2 = R^2  .
% \end{multline}
% After distributing and simplifying, we obtain a quadratic equation for $\overline{EI}$:
\begin{multline}
    \overline{EI}^2 + (2x_E\cos{\psi_ E})\overline{EI} \\  + (x_E^2 + \rho^2 - 2x_E\rho\sin{\psi_ E} - R^2) = 0 .
\end{multline}

Let $a=1$, $b=2x_E\cos{\psi_ E}$, and $c=x_E^2 + \rho^2 - 2x_E\rho\sin{\psi_ E} - R^2$, and substitute into the standard quadratic formula to obtain the equation for $\overline{EI}$ as: 
\begin{multline}
\label{eq:EI}
    \overline{EI} = -x_E\cos{\psi_ E} \pm  \\ \sqrt{-x_E^2\sin{\psi_E} + 2x_E\rho\sin{\psi_ E} - \rho^2 + R^2} .
\end{multline}

We now have an equation for $\overline{EI}$, from which only the positive solution represents capture when time is positive, and the negative case represents the condition where $E$ moves backwards rather than forwards. 
This distance is used to obtain the time travelled by $E$ before capture by substituting the above expression into~\eqref{eq:time_eq}.
%This can now be used to substitute $d$ in \eqref{eq:time_eq}.
Now, to calculate the distance travelled by $P$, we must first find the final location of $P$.
By assuming $\theta_F$ is known, we can set the final location angle of $P$ as simply an offset of $\phi$ from $\theta_F$, as shown in Fig.~\ref{fig:tan_geo}.
Related to this angle is the segment $\overline{PF}$, or $c$ as shown in Fig. \ref{fig:tan_geo}, as $c$ represents the chord of $\phi$.
We calculate $c$ using the Pythagorean theorem noting that $\overline{IF}=\overline{EF}-\overline{EI}$.
We obtain $c=\sqrt{\overline{IF}^2 + \overline{IP}^2}$
Now, we can simply solve for the offset angle, $\phi$, by using the chord formula: 
\begin{equation}
    \label{eq:phi_tan}
    \phi = 2\arcsin{\frac{c}{2R}}.
\end{equation}
The final angle of $P$ is described as: 
\begin{equation}
    \label{eq:theta_pf_tan}
    \theta_{Pf}=\theta_F - 2\arcsin{\frac{c}{2R}}.
\end{equation}
We calculate the initial position of $P$ for a TAC condition as follows: 
\begin{equation}
\label{eq:TC_pc}
\theta_{P0} = \theta_F -\frac{\gamma\overline{EI}}{R} - 2\arcsin{\frac{c}{2R}}.
\end{equation}

\begin{figure}[ht]
    \centering
    \includegraphics[width=0.8\linewidth]{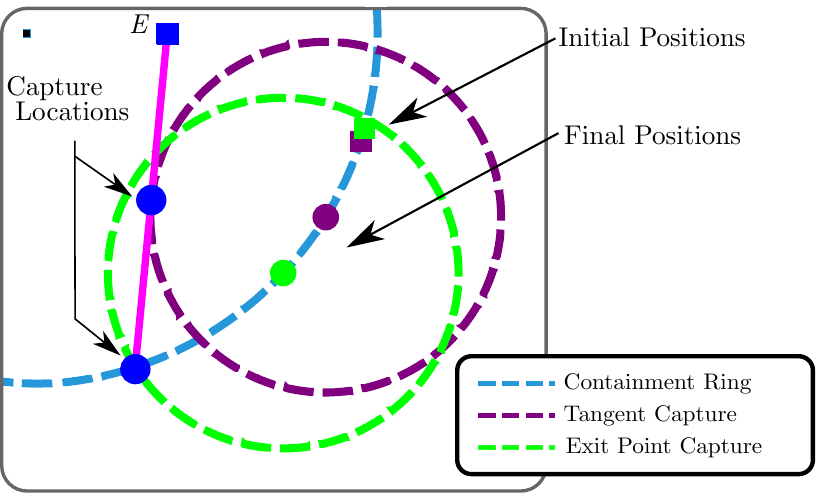}
    \caption{Comparison of TAC and EXC Conditions for Parameters: $\rho = 0.5$, $\gamma=0.5$, $\psi_E=1.6\pi$}
    \label{fig:comp_tan_vs_epc}
\end{figure}

\emph{Analysis:} After analyzing the TAC condition on different scenarios, it was observed that for all scenarios, the TAC condition resulted in a shorter range in initial conditions compared to the EXC condition shown before.
This is shown in Fig.~\ref{fig:comp_tan_vs_epc} where the square markers represent the initial positions for each of the capture conditions given.
Note the green square marker representing the EXC is placed further back along the ring compared to the purple square marker representing the TAC condition.
However, the final positions where capture occurs give us two limiting cases between which another capture condition may be found.

\subsubsection{Touch-and-Go Capture} The Touch-and-Go capture (TGC) condition occurs when two conditions are met: $R_{PE}=\rho$, the basic condition that assures $E$ has been captured, and $\dot{R}_{PE}=0$, the distance rate is zero. We can note from Fig. \ref{fig:epc_distance} that the blue curve has a local minima that represents where $\dot{R}_{PE}=0$ as shown in the orange curve of the same figure. The TGC condition captures the evader at the exact moment where the distance between $E$ and $P$ is at its minimum. The final position of $P$ is found numerically using the Bisection search method by setting the final locations of the TAC and EXC conditions as upper and lower boundaries and looking for the final position of $P$ that meets the two conditions stated. Once this final location is found, the position of $P$ is backtracked to find its initial location along the ring. 
\begin{figure}[ht]
    \centering
    \includegraphics[width=0.8\linewidth]{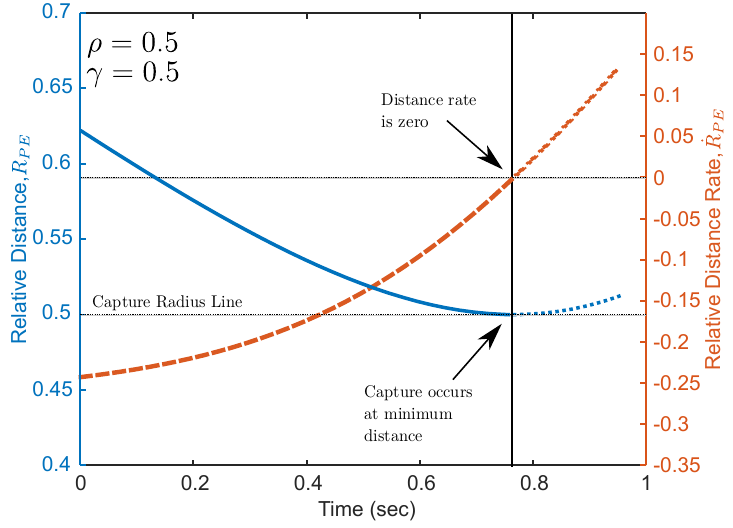}
    \caption{TGC Relative Distance and Relative Distance Rate }
    \label{fig:TGC_rdot plot}
\end{figure}

\begin{figure}
    \centering
    \includegraphics[width=0.45\linewidth]{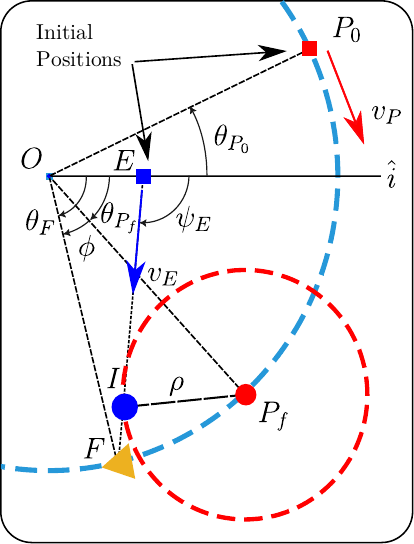}
    \caption{Touch-and-Go Capture Geometry}
    \label{fig:tng_geo}
\end{figure}

\begin{figure}[ht]
    \centering
    \includegraphics[width=0.8\linewidth]{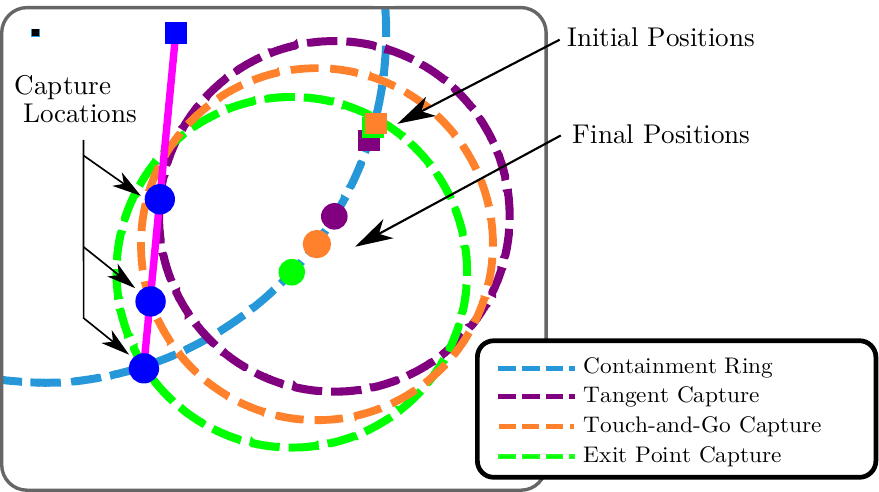}
    \caption{Comparison of TAC, EXC and TGC Conditions for Parameters: $\rho = 0.5$, $\gamma=0.5$, $\psi_E=1.6\pi$}
    \label{fig:comp_all_capture}
\end{figure}

\emph{Analysis:} After analyzing the TGC it was found that for a specific set of speed ratio, and capture radius size, there exists an interval of heading angles of E where the TGC conditions covers the most range (in terms of the pursuer's initial position) when compared to EXC and TAC. This is shown in Fig. \ref{fig:comp_all_capture} where TAC, EXC and TGC conditions are compared. Note the initial location of the TGC condition (orange square marker) is located further away along the ring compared to the other two conditions. Additionally, as the heading angle reaches the extremes of such interval, the TGC slowly transitions into a EXC condition. 
\subsection{Worst-case Capture}
\begin{figure}[ht]
    \centering
    \includegraphics[width=0.7\linewidth]{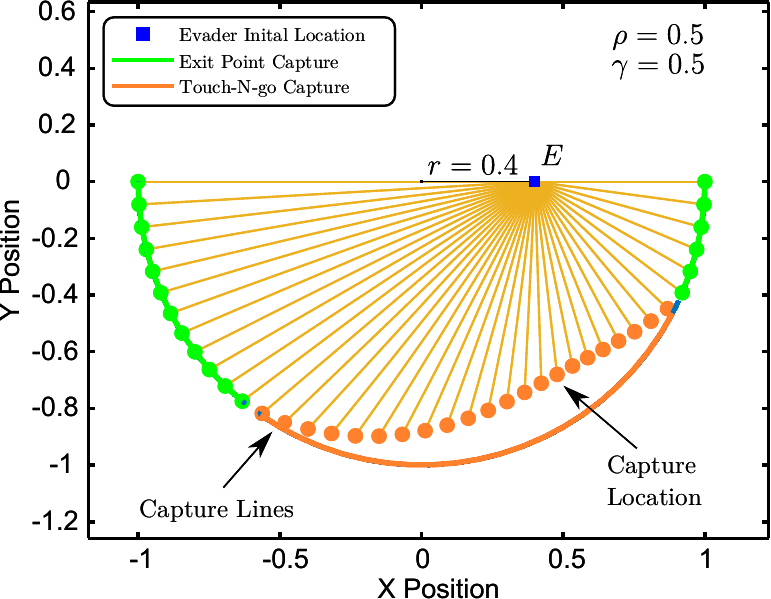}
    \caption{Worst-case Capture for Different Headings of $E$}
    \label{fig:opt_cap_ring}
\end{figure}

Worst-case capture for this problem is defined as the capture condition that covers the most angular range along the containment ring. This is the condition that allows the pursuer to be furthest back along the ring and still be able to contain $E$. This creates the most limiting scenario for $P$. Different types of capture conditions have been described and analyzed to define such boundaries. It was found that for a given set of fixed parameters $\rho$ and $\gamma$, the worst-case capture condition type depends on $\psi_E$, having EXC as worst-case when $\psi_E$ is close to $0$ or $\pi$, and TGC as worst-case for the interval in between. An example is shown in Fig.~\ref{fig:opt_cap_ring} where the set parameters are $\rho=0.5$ and $\gamma=0.5$, and $E$ is initially located $0.4$ along the $\hat{i}$ axis as shown by the blue square marker. The circle markers denote the capture locations for escape headings where $\pi \le \psi_E \le 2\pi$ with $P$ capturing in a CW motion. The green markers show EXC conditions while the orange markers show TGC conditions.   

\section{Parametric Study}
\label{sec:parametric_study}

\begin{figure}[ht]
    \centering
    \includegraphics[width=0.7\linewidth]{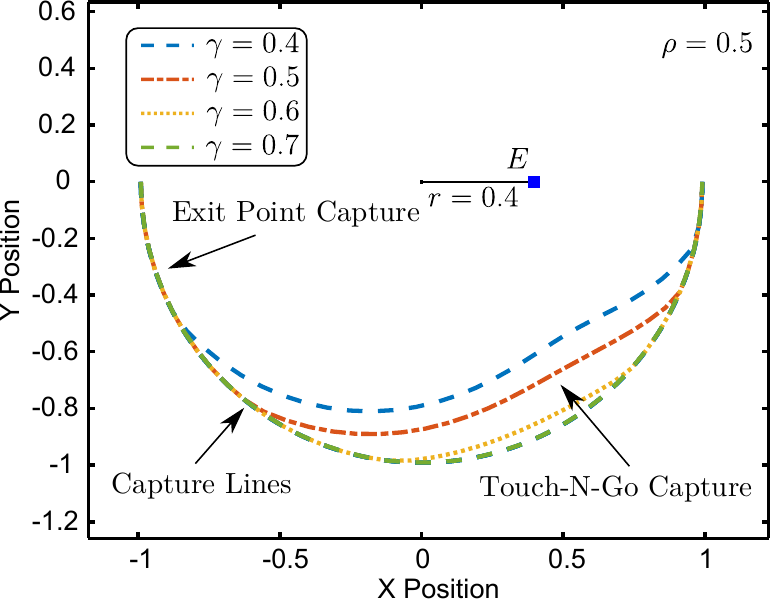}
    \caption{Capture Locations Varying Speed Ratio ($\rho=0.5$)}
    \label{fig:param study speed ratio}
\end{figure}

\begin{figure}[t]
    \centering
    \includegraphics[width=0.7\linewidth]{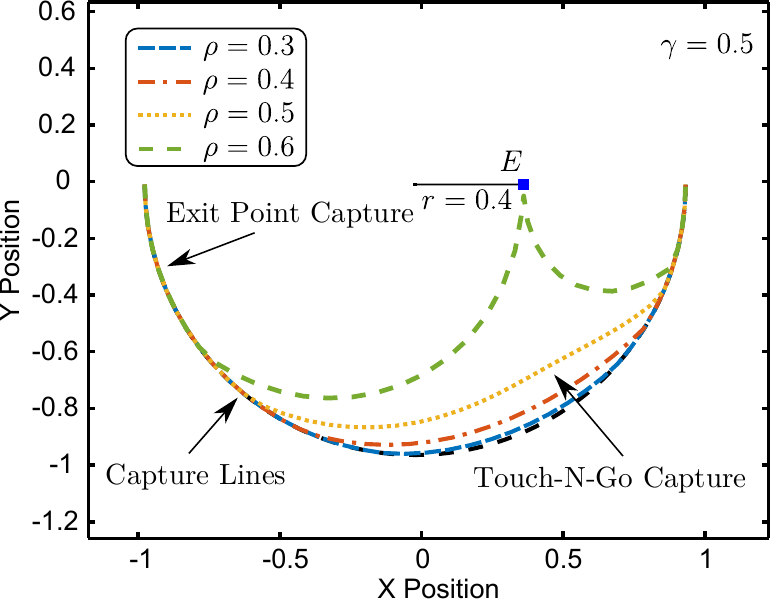}
    \caption{Capture Locations Varying Capture Radius ($\gamma=0.5$)}
    \label{fig:param study cap radius}
\end{figure}

\begin{figure}[hb]
    \centering
    \includegraphics[width=0.7\linewidth]{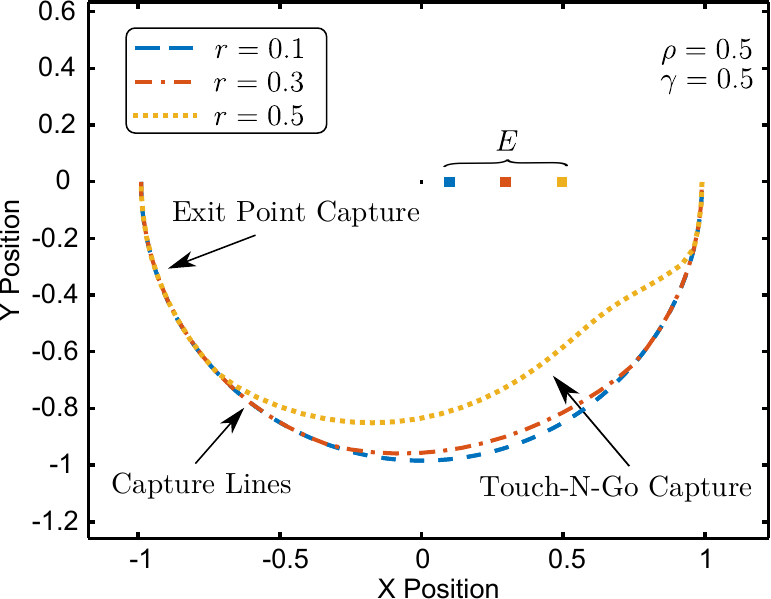}
    \caption{Capture Locations Varying Evaders Initial Location ($\rho=\gamma=0.5$)}
    \label{fig:param study re0}
\end{figure}

A parametric study is performed to investigate the effects of varying the speed ratio, capture radius, and initial location of $E$. The following figures depict the location of capture along the bottom half of the containment region. Figure \ref{fig:param study speed ratio} shows the effects of varying the speed ratio, $\gamma$ while keeping the capture radius and initial location of $E$ as $\rho=0.5$ and $r=0.4$, respectively. Note that as $\gamma$ increases, the capture location occurs nearer to the edge of the circle.    
Similarly, Fig. \ref{fig:param study cap radius} shows the effects of varying the capture radius, $\rho$ while keeping the rest of the parameters as: $\gamma=0.5$ and $r=0.4$. Note that as $\rho$ decreases, the capture location occurs nearer to the edge of the circle region, until point capture is achieved and all captures are of the EXC condition.   
Finally,  Fig. \ref{fig:param study re0} shows the effects of varying initial location of $E$ along the $\hat{i}$-axis while keeping the set parameters as: $\rho=0.5$ and $\gamma=0.5$. Note that as the location of $E$ translates to the center of the circle, the capture locations occur closer to the edge of the circle region until all captures are of the EXC type. 

\section{Reachability Analysis}
\label{sec:reachability_analysis}

\begin{figure}[ht]
    \centering
    \includegraphics[width=0.7\linewidth]{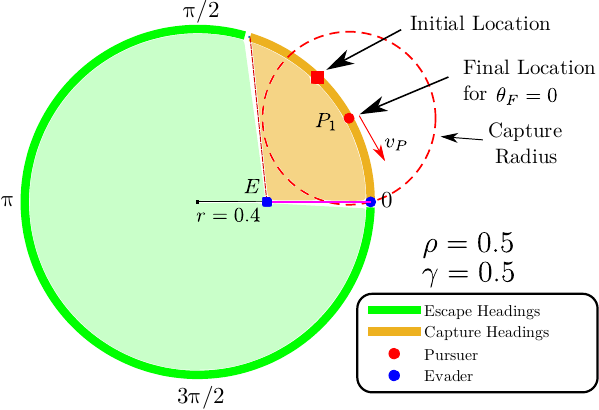}
    \caption{Reachability of One Pursuer}
    \label{fig:escape_1P}
\end{figure}

A reachability analysis is performed to study viable escape headings for $E$, i.e. the escape and capture regions of the problem. Figure \ref{fig:escape_1P} shows the basic scenario with a single evader and a single pursuer. $P$'s initial location is fixed as shown and is constrained to move CW. The green region represents the viable escape headings while the orange headings represents the capture headings. 

\begin{figure}[ht]
    \centering
    \includegraphics[width=1\linewidth]{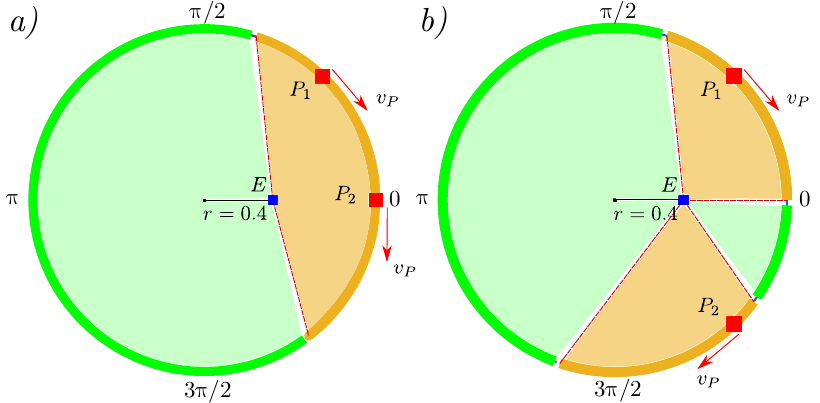}
    \caption{Reachability of Two Pursuers  \emph{a)} Consecutively Located, and  \emph{b)} Evenly Located}
    \label{fig:reach 2circles CW}
\end{figure}

Figure \ref{fig:reach 2circles CW} shows the case for a two-pursuer scenario where the pursuers move in CW motion and are \emph{a)} consecutively located along the ring and \emph{b)} evenly located with respect to $0$ rad along the ring.

\begin{figure}[ht]
    \centering
    \includegraphics[width=1\linewidth]{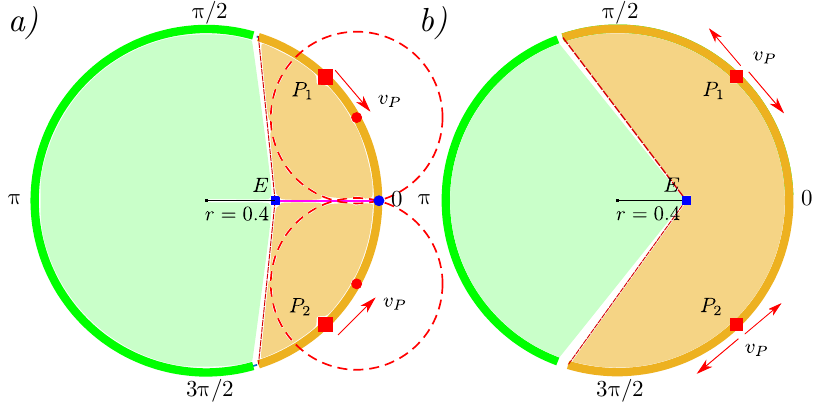}
    \caption{Reachability of Two Pursuers Evenly Located \emph{a)} Moving in Opposite Directions, and \emph{b)} Moving in Favorable Directions}
    \label{fig:reach 2circles CW CCW}
\end{figure}

Similarly, Fig. \ref{fig:reach 2circles CW CCW} shows the case for a two-pursuer scenario where the pursuers are evenly located with respect to $0$ rad and move in \emph{a)} a fixed CW and CCW motion and \emph{b)} the most favorable direction given $E$'s aim point angle, $\theta_F$, is known, that is, $P$ can decide to move CW or CCW depending on which will correspond to a smaller angular range between $P$'s initial position and $E$'s heading.  

\section{Conclusion} 
\label{sec:conclusion}

The problem of a single evader attempting to escape from a containment circular region and a pursuer attempting to capture the evader was considered. This scenario resembles the classical Lady in the Lake problem, but analyzes the case where the evader is faster than the pursuer, but the pursuer has a nonzero capture radius. Different capture conditions were discussed to define the most limiting case for capture to occur. A parametric study was performed to investigate the effects of speed ratio, capture radius size, and initial location of the evader. Finally, a reachability analysis was performed to find viable escape headings for the evader and illustrate escape and capture regions. As expected, a single pursuer is unable to guarantee capture for most headings, thus, the pursuer must rely on a team to extend the capture regions. This analysis represents foundational work for future extensions of the problem such as the incorporation of multi-pursuers, an evader with more complex motion, and malicious agents that aid the evader to escape.

\bibliographystyle{ieeetr}
\bibliography{references}

\end{document}